\documentclass[letterpaper, 10 pt, conference]{ieeeconf}  

\IEEEoverridecommandlockouts

\usepackage{graphicx} 
\usepackage{amsmath} 
\usepackage{amssymb}  
\usepackage{algorithm}
\usepackage[noend]{algpseudocode}
\usepackage{xcolor}
\newcommand{\fix}[1]{{\color{red}#1}}
\renewcommand{\fix}[1]{#1}
\title{\LARGE \bf
Iterative Switching Time Optimization for Mixed-integer Optimal Control Problems 
}

\author{Ramin Abbasi-Esfeden$^{1}$, Wim Van Roy$^{1}$, and Jan Swevers$^{2}$%
\thanks{$^{1}$ Ramin Abbasi-Esfeden and Wim Van Roy are with MECO research team, Dept ME, KU Leuven,  Leuven, Belgium. {\tt\small ramin.abbasiesfeden@kuleuven.be, wim.vanroy@kuleuven.be}}
\thanks{$^{2}$ Jan Swevers is with MECO research team, Dept ME, KU Leuven and Flanders Make@KU Leuven, Leuven, Belgium. {\tt\small jan.swevers@kuleuven.be}}%
}

\newcommand{\floor}[1]{\lfloor#1\rfloor}

\newcommand{\Sc}{\mathcal{S}}
\newcommand{\Wc}{\mathcal{W}}

\begin{document}

\maketitle
\thispagestyle{empty}
\pagestyle{empty}

\begin{abstract}
This paper proposes an iterative method to solve Mixed-Integer Optimal Control Problems arising from systems with switched dynamics. The so-called relaxed problem plays a central role within this context. Through a numerical example, it is shown why relying on the relaxed problem can lead the solution astray. As an alternative, an iterative Switching Time Optimization method is proposed. The method consists of two components that iteratively interact: a Switching Time Optimization (STO) problem and a sequence optimization. Each component is explained in detail, and the numerical example is resolved, the results of which shows the efficiency of the proposed algorithm. Finally, the advantages and disadvantages of the method are discussed and future lines of research are sketched.
\end{abstract}

\section{INTRODUCTION}

Switched dynamical systems show a large flexibility in modeling a wide range of real-world applications. However, due to the discrete nature of the switched dynamics, the optimal control of such a system is difficult. This paper proposes an approach to solve problems of the following form
\begin{subequations}
\begin{align}
\min_{x, u, v} &\int_{t_0}^{t_f} L\big(x(t), u(t), v(t),\fix{t}\big) dt + M\big(x(t_f)\big) \\
s.t. & \\
&\dot{x}(t) = f\big(x(t), u(t), v(t), \fix{t}\big),\\
&h\big(x(t), u(t), v(t), \fix{t}\big) \leq 0, \\
&x(t_0) = x_0, \\
& u(t) \in \{u_0, u_1, \ldots, \fix{u_{m}}\}, \;\fix{\forall t \in [t_0, t_f],}
\end{align}
\end{subequations}
where $x(t) \fix{\in \mathbb{R}^{n_x}}$ is the system state, and $u(t)\fix{\in\mathbb{R}^{n_u}}$ and $v(t) \fix{\in \mathbb{R}^{n_v}}$ are the inputs. The switching nature of the systems is captured in the variable $u(t)$, which can only take values from a finite discrete set of options. This is known as a Mixed-Integer Optimal Control Problem (MIOCP). We are only interested in efficient algorithms that can solve such a problem in real-time. For a complete survey of different approaches, see \cite{Zhu2015}.

A direct way to tackle this problem is to discretize the differential equations and solve this problem as Mixed-Integer Nonlinear Programming (MINLP). However, solving MINLPs, in general, can be difficult, and approximate solutions are sought as in \cite{Sager2012} and \cite{ Kirches2011}, where the approximation is done using rounding schemes based on the relaxed solution. The relaxed solution is an important modification of the original MIOCP, in which the discrete variable $u(t)$ is allowed to assume continuous values.

A different approach is to recognize two subproblems within the original one: (1) finding a switching sequence and (2) when a switching sequence is given, optimizing the switching instances; the rest can be cast as continuous nonlinear programming, for which efficient numerical algorithms exist. We refer to the second problem as {\it Switching Time Optimization (STO)} and call the first problem {\it Sequence Optimization}. 

Define a {\it sequence} as the sequence $\mathcal{S} = \{s_0, s_1, \ldots, s_{ns}\}$ of operation \fix{stages} (or modes) of the system ordered by their happening in time. Each operation \fix{stage} corresponds to a fixed value of discrete inputs. To specify the switching instance between two operation \fix{stages}, we define the set $\Wc = \{w_0, w_1, \ldots, w_{ns}\}$ based on $w_i$ defined as the time the system dwells in the \fix{stage} $s_i$.

STO problem is to calculate $\Wc$ given $\Sc$. Xu and Antsaklis \cite{Xu2002} use derivative information of the cost function to calculate the switching times, and similar treatments are found in \cite{Caldwell2011} and \cite{Egerstedt2003}. Lee et al. \cite{Lee1999} use a time transformation, also known as the control parametrization enhancing technique (CPET), to calculate $\Wc$. This is the method that is applied in this paper to solve STO. 

For STO to provide a satisfactory solution to the original MIOCP, an optimal $\Sc$ is to be proposed. Lee et al. \cite{Lee1999} suggest starting from an initial $\Sc$ and expanding the $\Sc$ until no decrease in the cost function is measured. Wardi and Egerstedt \cite{Wardi2012} start with an initial $\Sc$ and using the derivative information of the cost function, insert required \fix{stages} into $\Sc$. The insertion method is also found in \cite{Axelsson2008} and \cite{Gonzalez2010}.

Another family of solutions is based on the relaxed solution. One can show that a continuous solution can be adequately projected into a discrete solution \cite{Sager2012}\cite{Bengea2005}. This is used by Bengea and DeCarlo \cite{Bengea2005} to provide a sequence based on the relaxed solution. Gerdts and Sager \cite{Gerdts2012} use similar reasoning to find $\Sc$ (and $\Wc$) using a time transformation. Vasudevan et al. \cite{Vasudevan2013a}\cite{Vasudevan2013b} in a two parts paper present an algorithm based on the relaxed solution. Wei et al. \cite{Wei2007} in the same spirit as \cite{Vasudevan2013b} builds upon the relaxed solution.

In a more sophisticated way, Sager and Zeile \cite{Sager2021} use Mixed-Integer linear programming, also known as the combinatorial integral approximation (CIA), to project the relaxed solution into binary values. Based on this method, an open-source tool known as Pycombina \cite{Burger2020} is written to solve MIOCP problems. 

However, there are shortcomings with all the aforementioned methods. Expanding the initial $\Sc$ in \cite{Lee1999} assumes that the only constraints on $\Wc$ values are simple non-negativity constraints. However, there might be a nonzero lower bound on values of $\Wc$. Based on the results \cite{Wardi2012} provides, inserting \fix{stages} into $\Sc$ leads to frequent switching and cannot take constraints on $\Wc$ into account.

Other approaches rely heavily on the relaxed solution. Through a numerical example, we show that the the usage of the relaxed solution is justified only if chattering is not prohibited. If there exist limiting constraints on the values of $\Wc$, the methods that depend on the relaxed solution can lead to unsatisfactory results. 

In this paper, to find the optimum sequence, we start from an initial sequence $\Sc$ that is neither required to be optimal nor based on the relaxed solution, rather it only needs to {\it include} the optimal sequence as a set. The extra (i.e., suboptimal) \fix{stages} that are not part of the optimal sequence are iteratively removed from the initial sequence. As an application of the proposed method, a modified version of a benchmark problem known as Double Tank \cite{Vasudevan2013b} is solved.

\section{PROBLEM FORMULATION}

\subsection{The Double Tank problem}
\begin{figure}[thpb]
\centering
\makebox{\parbox{3in}{\centering
\includegraphics[width=1.3in]{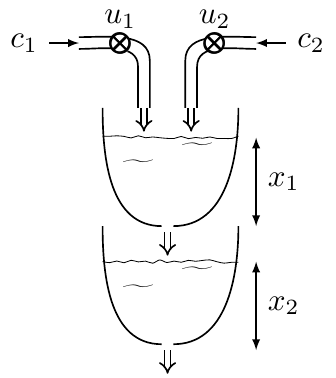}
}}
\caption{\fix{The Double Tank problem. The level of the second tank ($x_2$) is controlled using the valves ($u_1$, $u_2$) and the flow ($c_1$, $c_2$) of two pipes.}}
\label{fig:doubletank}
\end{figure}

\fix{The Double Tank problem is shown in Fig. \ref{fig:doubletank}. Instead of two modes for the same pipe, as in the original Double Tank \cite{Vasudevan2013b}, }there are two pipes, one of which has a variable parameter for its flow. The level of the second tank is expected to follow a sine function as a reference, and there is also an added cost, albeit small, for using each pipe. \fix{Each pipe has a valve that can open or close it, and both pipes can be used simultaneously. The problem is formulated as }

\begin{subequations}
\begin{align}
\min_{x,u, c} \int_0^{t_f} \Big[ \alpha &(x_2(t) - r(t))^2 + \sum_{i=1}^2 \beta_{i} u_i(t) c_i(t) \Big] dt \\
s.t. \quad \dot{x}_1(t) &= \sum_{i=1}^2 u_i(t) c_i(t) - \sqrt{x_1(t)},\\
 \quad \dot{x}_2(t) &= \sqrt{x_1(t)} - \sqrt{x_2(t)}, \\
&  \fix{c_1(t) = \gamma},\\
0 & \leq c_2(t) \leq \gamma, \\
u_i(t) &\in \{0,1\}\quad \forall i \in \{1, 2\}, \label{eq:bin}
\end{align}
\end{subequations}
with a reference trajectory $r(t) = 2 + 0.5 \sin(t)$, $x(0) = (2,2.5)$ and constants
\begin{align*}
t_f = 10 \;,\; \alpha = 100 \;,\; \beta_1 = 1 \;,\; \beta_2 = 1.1 \;,\; \gamma = 10.
\end{align*}

\subsection{The relaxed solution}
\begin{figure}[thpb]
\centering
\makebox{\parbox{3in}{
\includegraphics[width=3in]{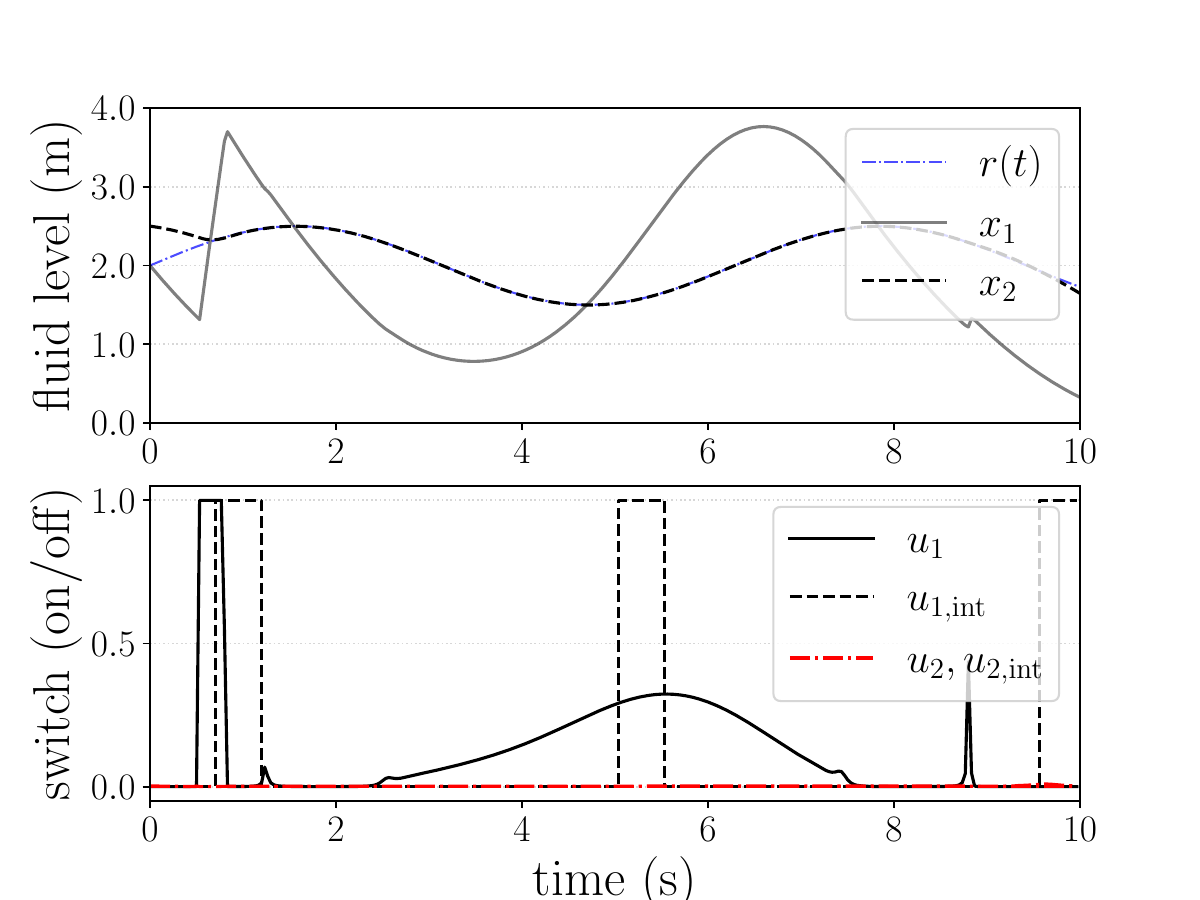}
}}
\caption{The relaxed solution to the modified Double Tank problem. The top plot shows the fluid level of the second tank following the reference exactly. The bottom plot shows the relaxed values for switches, along with their projected values.}
\label{fig:rel}
\end{figure}

To see the shortcoming of the approaches that depend on the relaxed solution, the relaxed solution is calculated, and it is projected into an integer solution through CIA formulation.  The relaxed solution and its integer projection are shown in Fig. \ref{fig:rel}. The solution is only using $u_1$, which is the cheaper option, to make $x_2$ follow the sine reference, and it does so successfully. 

\fix{The relaxed problem is discretized in a multiple shooting scheme for 300 nodes, and for simplicity, the Explicit Euler method is used as the integrator. The resulting nonlinear programming is then solved} using the interior-point solver IPOPT \cite{Wachter2005} with a runtime of $0.5\mathrm{(s)}$. It is then projected into an integer solution using Pycombina BnB solver in $0.8\mathrm{(s)}$. To prevent chattering, a minimum uptime of $0.5\mathrm{(s)}$ is specified.

\begin{figure}[thpb]
\centering
\makebox{\parbox{3in}{\centering
\includegraphics[width=3in]{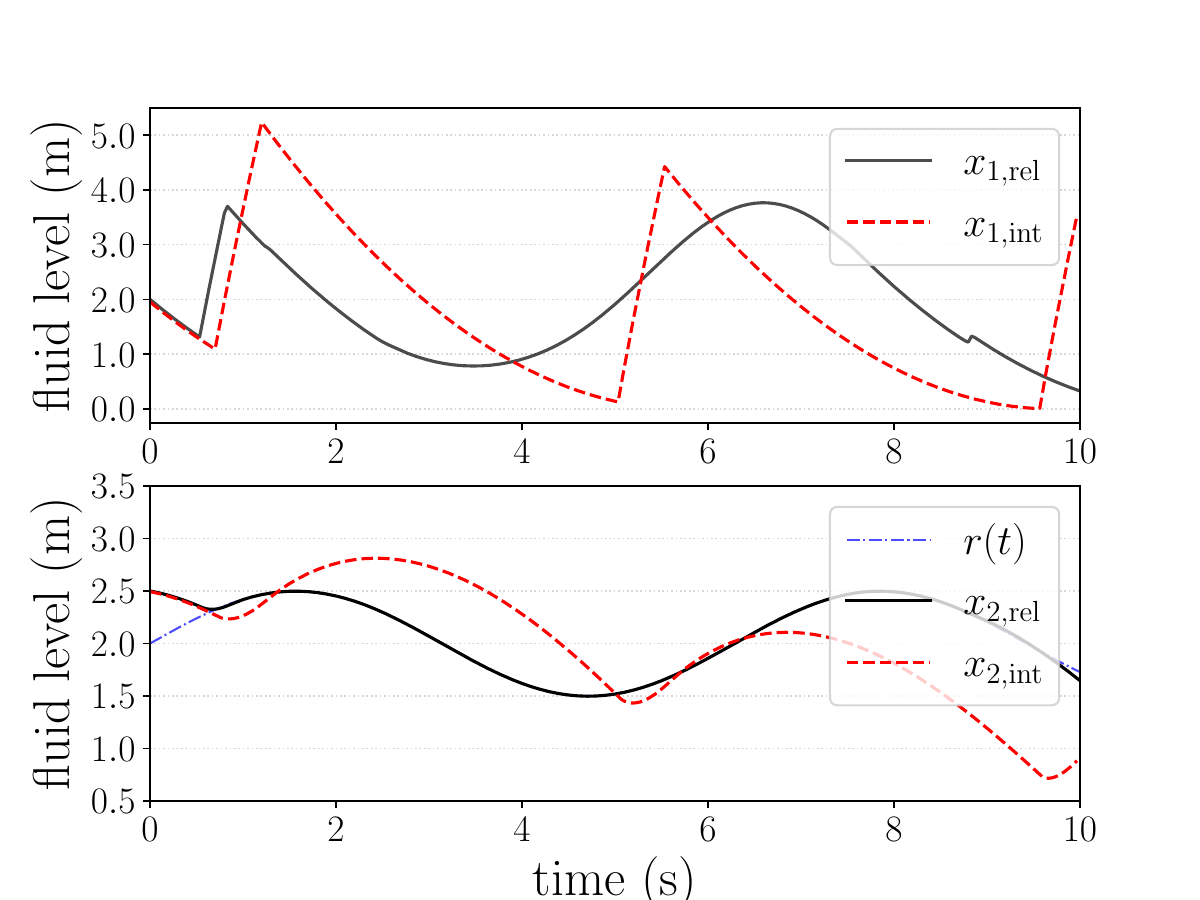}
}}
\caption{The fluid levels for the relaxed solution and its integer counterpart. The integer solution fails to follow the reference as there is not enough chattering allowed.}
\label{fig:relbin}
\end{figure}

Note that the relaxed solution is using the pipe that has a fixed flow $\gamma$. However, since it relaxes the switching variable $u_1(t)$, it can provide the necessary flow. The relaxation is giving $u_1(t)$ a flexibility that it does not possess, and this is bound to cause issues if no chattering is allowed. This issue becomes apparent if the outcome of Pycombina is applied and simulated. In Fig. \ref{fig:relbin} we see that the integer solution is not able to stir the $x_2$ to follow the reference. 

This is not satisfactory as we see that there is an obvious answer of $u_1(t) = 0, u_2(t) = 1$ to this problem, which would outperform the projected answer. The flexibility of the relaxed solution along with the greediness of the optimizer causes the integer solution to lose sight of this answer. Moreover, after fixing the integer solution, there remains no free variable to optimize further, as $u_2(t)$ is identically zero.



%

\section{ITERATIVE SWITCHING TIME OPTIMIZATION}
\subsection{Switching Time Optimization} 
\begin{figure}[thpb]
\centering
\makebox{\parbox{3in}{\centering
\includegraphics[width=3in]{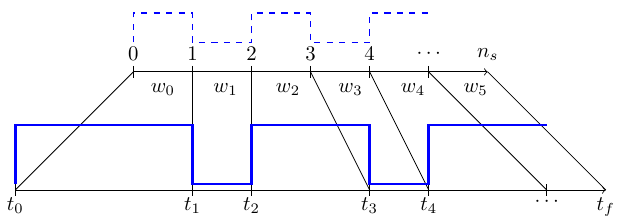}
}}
\caption{The time transformation from $\tau$ to $t$ through $w(\tau)$.}
\label{fig:stott}
\end{figure}

To solve the STO problem, we follow the time transformation approach of Lee et al. \cite{Lee1999}. Based on $\Wc$, we define the following function on the interval $[0, n_s]$ 
\begin{align}
 w(\tau) = w_i \;;\quad i \leq \tau < i+1,\; \forall i \in \{0,1, \ldots,n_s-1\}.
\end{align}
A time transformation from the interval $[0, n_s]$ onto $[0, t_f]$ can then be defined as
\begin{align}
\frac{dt}{d\tau} = w
\end{align}
The time transformation is visualized in Fig. \ref{fig:stott}. 

Using this time transformation, we can define the following function
\begin{align}
\begin{split}
t(\tau') &= \int_0^{\tau'} w d\tau, \\
&= \int_0^1 w_0d\tau + \int_1^2 w_1d\tau + \cdots + \int_{\floor{\tau'}}^{\tau'} w_{\floor{\tau'}} d\tau, \\
&= w_0 + w_1 + \cdots + (\tau' - \floor{\tau'}) w_{\floor{\tau'}},
\end{split}
\end{align}
and its inverse as
\begin{align}
\tau(t') = \min \{\tau | t(\tau) = t' \}.
\end{align}
This transformation helps to transfer general expressions such as
\begin{align}
\int_0^{t'} f(t)dt = \int_0^{\tau(t')} f(t(\tau)) w d\tau ,
\end{align}
which transforms the original problem into
\begin{subequations}
\begin{align}
\min_{x, w, c} \int_0^{n_s} \Big[ \alpha &(x_2(\tau) - r(\tau))^2 + \sum_{i=1}^2 \beta_{i}\bar{u}_i c_i(\tau) \Big] wd\tau \\
s.t. \quad \dfrac{dx_1(\tau)}{d\tau} &= w \Big[ \sum_{i=1}^2 \bar{u}_i c_i(\tau) - \sqrt{x_1(\tau)} \Big],\\
 \quad \dfrac{dx_2(\tau)}{d\tau} &= w\Big[ \sqrt{x_1(\tau)} - \sqrt{x_2(\tau)}\Big], \\
&  \fix{c_1(\tau) = \gamma},\\
0 & \leq c_2(\tau) \leq \gamma,\\
0 &\leq h(w(\tau)) \label{eq:hw},\\
t_f &= \int_0^{n_s} w(\tau) d\tau \label{eq:wsum}.
\end{align}
\end{subequations}

Note that the integer variables $u$ are now removed as a variable and enter the equations as constants $\bar{u}_i$ based on the \fix{stage} the system is in, which itself is determined by the value of $\tau$. Also, the final constraint makes sure that the $w$ is indeed a time transformation of the original system. Constraints on the switching times enter naturally into optimization by $h(w)$.

\subsection{Sequence optimization}
For sequence optimization, we start from a sequence $\Sc$ that is rich enough to include the optimal sequence $\Sc^*$, i.e. $\Sc^* \subset \Sc$. The main problem is to recognize the extra \fix{stages}, which is denoted by $\Sc^0 \subset \Sc$, and remove them from the sequence. This is to be done iteratively: following each removal, the STO is solved again for the new sequence, from which other candidates for removal are recognized. As the initial sequence is finite, the iteration comes to an end eventually.

The time transformation in STO formulation shows that \fix{stages} with zero $w$ are practically inert as they do not endure in time to have any effect on the system; they are an obvious candidate for redundant \fix{stages} and can be safely removed from the sequence. In other words, if $w_i = 0$, then $s_i \in \Sc^0$.

\begin{figure}[thpb]
\centering
\makebox{\parbox{3in}{\centering
\includegraphics[width=2.5in]{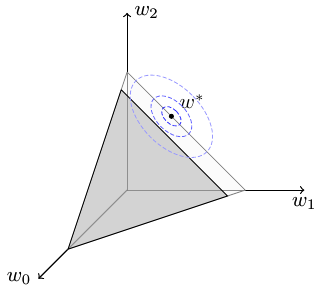}
}}
\caption{An example of a three \fix{stages} system with a lower bound constraint for $w_0$ with an optimum $w^* = (0, w_1^*, w_2^*)$. The grey polygon on the plane shows the feasible set for $(w_0, w_1, w_2)$.}
\label{fig:stoplane}
\end{figure}

If there were no constraints on the values of $\Wc$, we would expect the optimizer to push the extra \fix{stages} toward zero. However, there might be constraints that prevent some $w_i$ from becoming zero. An example of three \fix{stages} is shown in Fig. \ref{fig:stoplane}. Because of \eqref{eq:wsum}, the feasible values of $\Wc$ constitute a hyperplane. This hyperplane is restricted by \eqref{eq:hw} to non-negative values for $w_0, w_1, w_2$ and a lower bound on $w_0$. 

A workaround can be recognized in Fig. \ref{fig:stoplane}. At each iteration, any \fix{stage} that is not pushed against such constraints, by definition, is deemed optimal by the optimizer. Thus if there are any \fix{stages} to be removed, they are to be found within \fix{stages} with a nonzero Lagrange multiplier for the corresponding constraint. We will denote them by $\bar{\Sc} \subset \Sc$. 

Within $\bar{\Sc}$, we can distinguish the optimal \fix{stages} from the extra ones from the optimizer behavior if there were no constraints limiting its action. In that case, if the optimizer drives a \fix{stage} $s \in \bar{\Sc}$ to zero then $s \in \Sc^0$. If the \fix{stage} is not driven to zero, it could be either extra or optimal.

To turn the above discussions into an iterative scheme, slack variables $e_i$ are added to each \fix{stage} $s_i$ in the formulation of the limiting constraints. A corresponding slack cost $\phi_i \triangleq \fix{\frac{1}{2}a_ie_i^2}$ is then added to STO; this will help us to identify $\bar{\Sc}$. An additional cost $\psi_i \triangleq \fix{\frac{1}{2}b_i w_i^2}$ is also attached to STO; this is the required cost for driving \fix{stages} to zero.

We start by setting $a_i$ to a small nonzero value, and $b_i=0$. Next, STO is solved and any \fix{stage} with $w = 0$ is removed as they are in $\Sc^0$. This will help us to identify $\bar{\Sc}$ as any \fix{stage} with $e_i > 0$ is in $\bar{\Sc}$. If $\bar{\Sc}$ is empty, the optimal sequence is reached. Assuming that $\bar{\Sc}$ is not empty, denote the \fix{stage} with the largest slack $s_r$ as the candidate to be removed.

To check whether $s_r \in \Sc^0$ or not, we alternate the $\phi$ and $\psi$ costs by setting $a_r = 0$ and $b_r$ to a nonzero value. STO is solved again, and any \fix{stage} with $w = 0$ is removed. If $w_r = 0$, it was indeed in $\Sc^0$ and will be removed at this point, and the iteration moves to the next candidate. 

If a \fix{stage} is repeated as a candidate, it means that the initial values for $a_r$ and $b_r$ had not been enough to settle between constraints satisfaction or \fix{stage} removal. In such case we continue to alter values of $a_r$ and $b_r$, increasing them each time the same \fix{stage} is chosen as $s_r$. 


It is also possible to manually enforce either constraint satisfaction or \fix{stage} removal, but the alternation between costs for the undecided \fix{stages} ensures that not much is lost by keeping an extra \fix{stage} or removing a necessary \fix{stage}. An example of such alternation for $a_r$ and $b_r$ is shown in Table. \ref{tab:ab}.  \fix{The above discussion is summarized in Algorithm \ref{iSTO}.}

\fix{
Algorithm \ref{iSTO} will not terminate only if there is always a \fix{stage} $s_r$ for which either $w_r = 0$ or $e_r > 0$. This implies either an infinite number of \fix{stage} removal or an infinite number of alternations of $a_r$ and $b_r$. However, if the problem is feasible and the initial $\Sc$ is finite, both cases are impossible. In the former case, it will require infinite \fix{stage} removal from a finite sequence, which is impossible. In the latter, it requires alternating values of $a_r$ and $b_r$ to ever-increasing values. This indicates that there is a \fix{stage} that cannot be driven to zero and cannot be made to satisfy a constraint regardless of the cost, which means that the original problem was infeasible.
}

\begin{algorithm}
\caption{\fix{Iterative Switching Time Optimization}}\label{iSTO}
\begin{algorithmic}[1]
\State Receive an initial sequence $\Sc$. 
\State Formulate a STO problem based on $\Sc$.
\State Soften timing constraints on $w_i$ by introducing slack variable $e_i$ for all $i \in \{1,2,\cdots,ns\}$.
\State Add $\sum_{i}^{ns}(\phi_i + \psi_i)$ to the STO cost function.
\State Set $a_i$ and $b_i$ based on Table. \ref{tab:ab}, starting from the first column.
\While {TRUE}
\State Solve STO for the sequence $\Sc$.
\If {$e_i = 0$ for all $i \in \{1,2,\cdots,ns\}$}
\If {$w_i > 0$ for all $i \in \{1,2,\cdots,ns\}$}
\State End.
\EndIf
\Else
\State Let $e_r$ be the largest slack variable.
\State Alternate $a_r$ and $b_r$ values based on Table. \ref{tab:ab}.
\State Solve STO for the sequence $\Sc$.
\EndIf
\State Remove \fix{stages} with $w = 0$ from the sequence $\Sc$. 
\EndWhile
\end{algorithmic}
\end{algorithm}

\subsection{Solution}
The iteration of STO and sequence optimization allows one to start from an initial sequence that includes an optimal sequence and reach it. Thus if we start from a sequence with $n_s$ element, there are $2^{n_s} - 1$ (the number of its power set elements excluding the empty set) possible sequences derivable from it. To solve the Double Tank problem, the following initial sequence of discrete inputs is chosen
\begin{align}
\Sc_0 = \{
\begin{bmatrix}
1\\
1
\end{bmatrix},
\begin{bmatrix}
0\\
1
\end{bmatrix},
\begin{bmatrix}
1\\
0
\end{bmatrix},
\begin{bmatrix}
0\\
0
\end{bmatrix},
\begin{bmatrix}
1\\
1
\end{bmatrix},
\begin{bmatrix}
0\\
1
\end{bmatrix},
\begin{bmatrix}
1\\
0
\end{bmatrix}
\}.
\end{align}
\renewcommand{\arraystretch}{1.2}

\begin{figure}[thpb]
\centering
\makebox{\parbox{3in}{\centering
\includegraphics[width=3in]{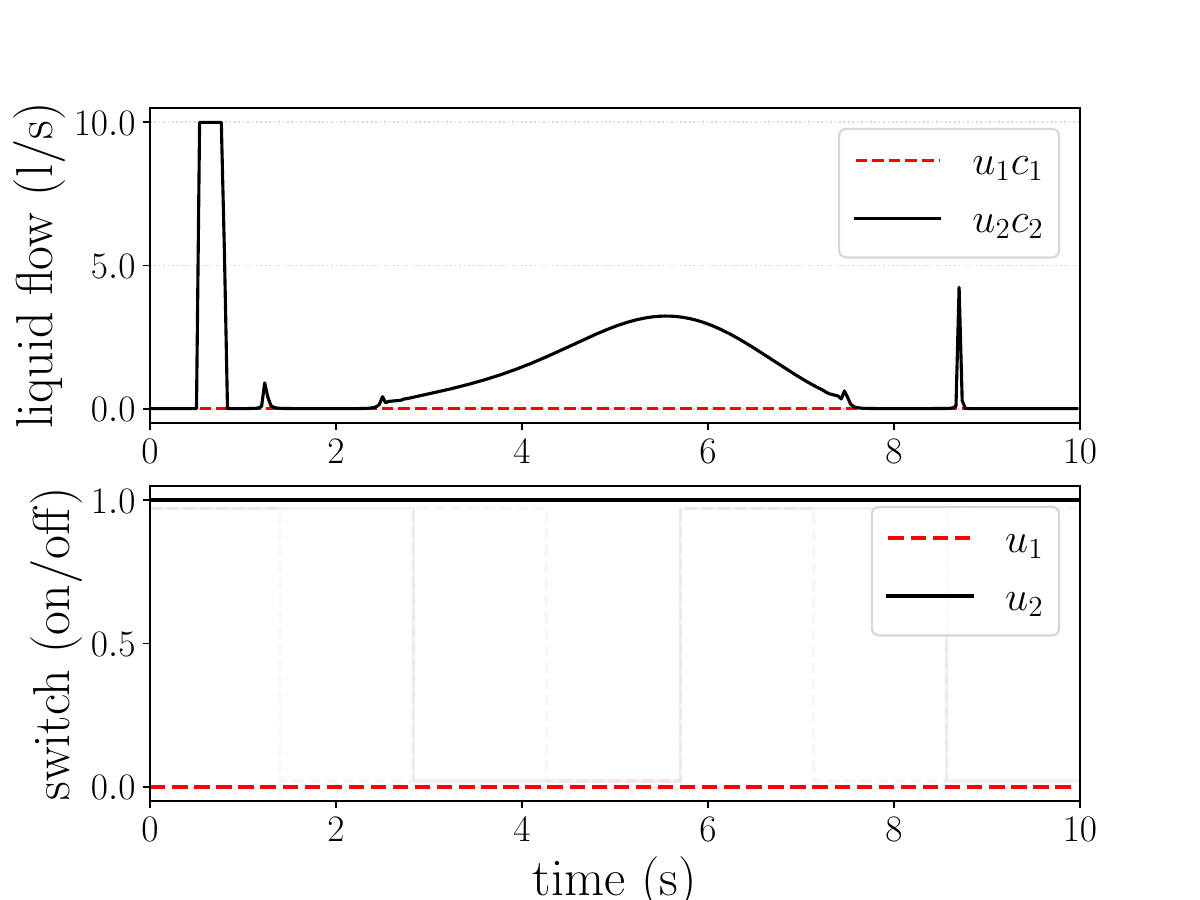}
}}
\caption{The iterative STO solution to the example problem. It produces the same flow as the relaxed solution, using $u_2$, however.}
\label{fig:sto}
\end{figure}

If we set a minimum uptime of $0.5\mathrm{(s)}$ and apply the iterative STO algorithm, after six iteration, the expected answer of $u_1 = 0, u_2 = 1$ is reached, as shown in Fig. \ref{fig:sto}. Six \fix{stages} were deemed redundant or infeasible by the optimizer and removed from the pool of \fix{stages}. 

\begin{table}[thpb]
\caption{An Example of a series of costs.}
\label{tab:ab}
\begin{center}
\begin{tabular}{|c|ccccc|}
\hline 
Repetition of $s_r$ & \fix{0} & 1& 2& 3 & 4\\
\hline
$a_r$ &$\fix{10^0}$& $0$ & $10^2$ & $0$ & $10^4$\\
$b_r$  &$\fix{0}$ & $10^1$ &$0$ & $10^3$ & $0$\\
\hline
\end{tabular}
\end{center}
\end{table}

Note that the peak value around $t = 0.5\mathrm{(s)}$ could have been produced by $u_1$, and it would have been cheaper as $u_1$ has a lower cost $\beta$. However, the minimum uptime made that answer infeasible, and the \fix{stage} was driven to zero and discarded. This can be confirmed if the minimum uptime constraint is lifted. 

\begin{figure}[thpb]
\centering
\makebox{\parbox{3in}{\centering
\includegraphics[width=3in]{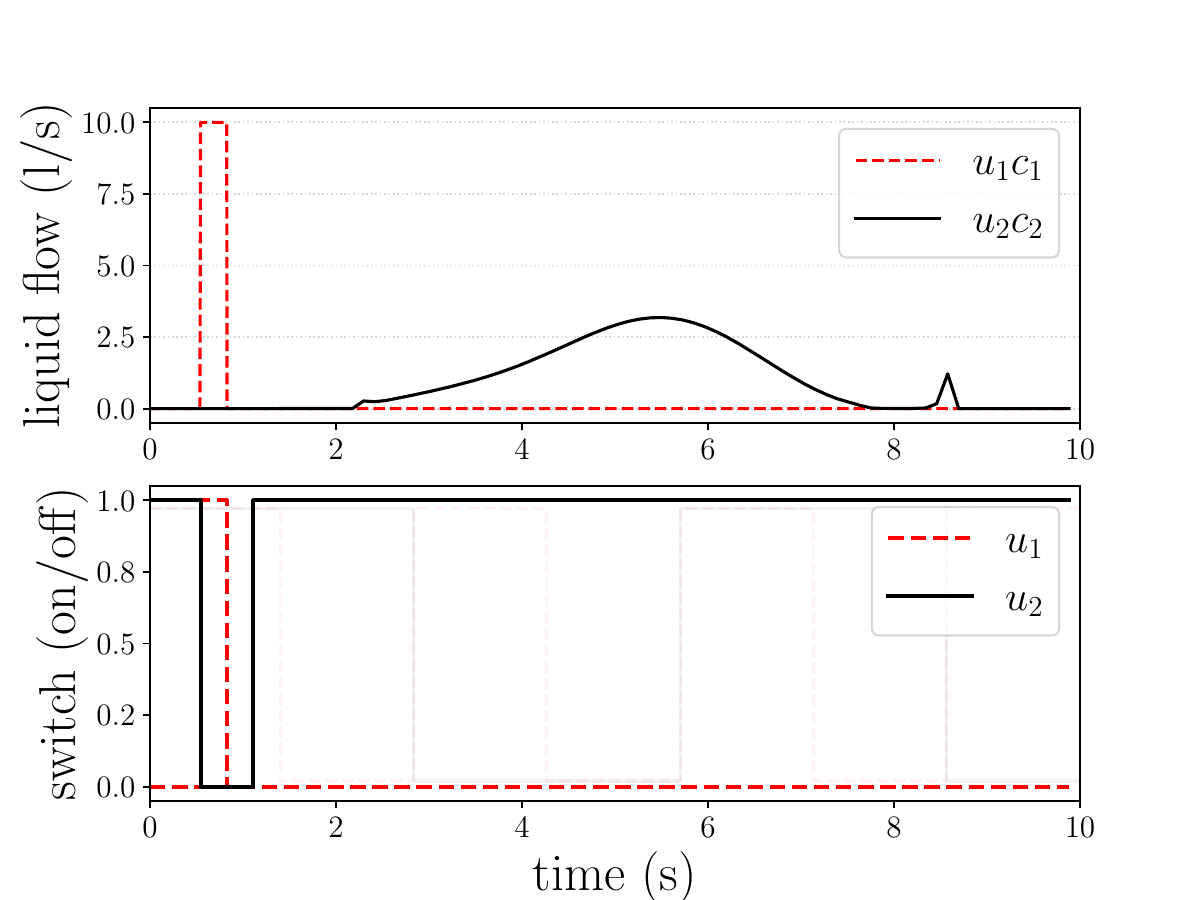}
}}
\caption{The iterative STO solution to the example problem without the minimum uptime constraint. The peak value is produced by $u_1$, which is the cheaper option.}
\label{fig:sto2}
\end{figure}

The results in Fig. \ref{fig:sto2} now shows the usage of the cheaper pipe $u_1$. This causes a cost reduction from $19.406$ to $18.702$. The relaxed solution cost is $18.239$ (the cost of the projected solution from Pycombina is over $300$ as it fails to follow the reference.) it takes $6$ iterations with the following runtime (total of $1.32\mathrm{(s)}$) to reach the solution
\begin{align}
\text{runtime (s)} = &
\begin{bmatrix}
0.53, 0.20, 0.11, 0.10, 0.26, 0.12
\end{bmatrix}.
\end{align}
Without the minimum uptime, this is reduced to $3$ iterations and the runtime of $[0.54, 0.13, 0.11]\mathrm{(s)}$. The first iteration takes longer as it starts from an initial guess, but later iterations are warm-started from the previous iterations.

The results show that iterative STO solves the problem in an efficient way without relying on the relaxed solution, thus avoiding the pitfalls discussed in the beginning. The efficiency of STO comes from the fact that it changes the mixed-integer problem into a continuous one in such a way that many constraints can enter the problem in a natural and simple way, e.g. $w > 0.5$ for minimum uptime.

There is no risk of frequent switching as the maximum number of switches is already limited by the initial sequence. There is also the possibility of applying conditional constraints such as {\it do not set $u_1 = 1$ unless $u_2 = 1$ or $u_3 = 1$ beforehand}, as they can be applied directly to the sequence to filter the prohibited combinations.

A common objection to such a method is the large number of possible combinations that are required for the initial sequence. In practice, this does not pose a big issue. First, the initial sequence only needs to {\it include} the sequence we seek, and the power set grows exponentially. Second, there does not have to be an increase in the number of variables. \fix{After discretization, STO fixes the values of $u(t)$ at each node and substitutes them with a small number of $w$, which reduces the number of variables in comparison to the relaxed problem.}

\section{FUTURE WORK}
The simplicity and efficiency of the proposed algorithm justify investigating it further. The efficiency of the algorithm and the quality of the answer relies on the initial sequence. A good sequence optimization cannot fix a bad initial sequence, so other sources for the initial sequence rather than simple brute enumeration needs to be investigated. 
\fix{One source can be the relaxed solution if the aforementioned risks are taken into account. Note that even though the relaxed solution is misusing the pipes, it is nonetheless correct in calculating the required input from the pipes for the liquid level to follow the reference; it is merely mistaken in allocating the right amount to each pipe. It is possible to sample the correct required input, derived from the relaxed solution, and among all the possible combinations of pipe usage that can produce that input, select the combination with the lowest cost. When considering switching cost or timing constraints, dynamic programming techniques can be applied to find the best combination. Therefore, it is plausible to use the relaxed solution to provide an initial sequence and optimize it further using iterative switching time optimization.}

\fix{The process of sequence optimization can also be optimized further. Note that if a suitable numerical optimization package is selected, STO can only lead to a decrease in the cost function. The same is true for removing stages with zero $w$, as they do not affect the system. In fact, if any switching cost is included, the stage removal will reduce the total cost. There is, however, a possibility that for the stages in $\bar{\Sc}$, which require alternation between $a_r$ and $b_r$, the algorithm leads to an increase in the cost because of removing a necessary stage or keeping a suboptimal one. As explained in the sequence optimization section, by limiting the increase in $a_r$ and $b_r$ at each alternation, it is possible to put an upper bound on this error. Alternatively, there can be a branch and bound scheme regarding the choice between constraint satisfaction and \fix{stage} removal by solving both options and comparing their cost to ensure the optimality of the sequence. Needless to say, the ability to insert a necessary \fix{stage} can bring further flexibility to the algorithm.}

\fix{Finally, the algorithm was applied to the Double Tank problem in an open loop fashion. The efficiency of the algorithm makes it a suitable option for employment in a closed loop Model Predictive Control. }

\addtolength{\textheight}{-13cm}   

\section{CONCLUSIONS}
In this paper, we proposed an iterative approach to solving Mixed-Integer Optimal Control Problems arising from systems with switched dynamics. Through a numerical example, it was shown why relying on the relaxed solution to provide an approximate projected answer might steer us toward unsatisfactory results. Two components of the method,  STO and sequence optimization were discussed in detail. The problem was solved using the proposed algorithm and the efficiency of the algorithm was shown from the optimum value of the cost function and the runtime required.
\section*{ACKNOWLEDGMENT}
This project has received funding from the European Union’s Horizon 2020 research and innovation program under the Marie Skłodowska-Curie agreement No.~953348

\end{document}